\documentclass{amsart}

\usepackage{amsmath,amssymb,amsthm,mathrsfs}
\usepackage[alphabetic]{amsrefs}
\usepackage[all]{xy}
\newdir{ >}{{}*!/-5pt/@{>}} %% cf xyguide exercise 14
\SelectTips{cm}{} %% cf xyguide page 9

\title{Chromatic Redshift}
\author{John Rognes}
\address{Department of Mathematics, University of Oslo, Norway}
\email{rognes@math.uio.no} \urladdr{http://folk.uio.no/rognes}
\date{January 28th 2014}

\newtheorem{theorem}{Theorem}[section]
\newtheorem{conjecture}[theorem]{Conjecture}

\newcommand{\bC}{\mathbb{C}}
\newcommand{\bF}{\mathbb{F}}
\newcommand{\bG}{\mathbb{G}}
\newcommand{\bH}{\mathbb{H}}
\newcommand{\bQ}{\mathbb{Q}}
\newcommand{\bW}{\mathbb{W}}
\newcommand{\bZ}{\mathbb{Z}}

\newcommand{\sA}{\mathscr{A}}
\newcommand{\sC}{\mathscr{C}}
\newcommand{\sE}{\mathscr{E}}
\newcommand{\sK}{\mathscr{K}}

\DeclareMathOperator{\Spec}{Spec}
\DeclareMathOperator{\Ext}{Ext}
\DeclareMathOperator{\Map}{Map}
\DeclareMathOperator{\Mod}{Mod}
\DeclareMathOperator{\THH}{THH}
\DeclareMathOperator*{\holim}{holim}
\newcommand{\longto}{\longrightarrow}
\renewcommand{\:}{\colon}

\DeclareMathOperator{\et}{\text{\'et}}
\DeclareMathOperator{\mot}{mot}

\begin{document}
\begin{abstract}
Notes for the author's MSRI lecture in January 2014.
\end{abstract}
\maketitle{}

\section{Introduction}

Consider a commutative ring $R$, with sum and product operations.
The category of representations of $R$ inherits a commutative rig
structure, given by direct sum and tensor product.  In other words,
the category $\Mod(R)$ of $R$-modules inherits a bipermutative
structure.  Continuing, one can consider the categorical representations
of $\Mod(R)$, and these in turn form a $2$-category $\Mod(\Mod(R))$,
with a ring-like structure.  Iterating, one can consider an $n$-category
of higher representations, for each $n\ge1$.

All of these constructions can take place within the limiting context of
structured ring spectra, or commutative $S$-algebras.  From the
category of (finite cell) modules over a commutative $S$-algebra $B$
we can distill a new commutative $S$-algebra, the algebraic
$K$-theory spectrum $K(B)$.  Continuing, one can form $K(K(B))$, etc.
When $B = HR$ is the Eilenberg--Mac\,Lane spectrum of an ordinary
ring, the $n$-fold algebraic $K$-theory $K^{(n)}(B)$ is extracted from
the $n$-category of higher representations, as above.
In this sense, $n$-fold iterated algebraic $K$-theory has something
to do with $n$-categories.

From this point of view it is surprising that $n$-fold iterated
algebraic $K$-theory also has something to do with formal group laws
of height~$n$, i.e., one-dimensional commutative formal group laws
$F$ in characteristic~$p$ where the series expansion $[p]_F(x)$ for
the multiplication-by-$p$ map starts with a unit times~$x^{p^n}$.
This is essentially a statement about the formal coproduct on
$K^{(n)}(B)^*(\bC P^\infty)$ that comes from the product on $\bC P^\infty$.
Hesselholt--Madsen asked about the chromatic filtration of iterated
topological cyclic homology in \cite{HM97}*{p.~61}, but could almost
as well have asked about the chromatic filtration of iterated algebraic
$K$-theory.

In a strong form, this connection implies that the algebraic $K$-theory of
a structured ring spectrum related to formal group laws of height~$n$
will be related to formal group laws of height~$n+1$.  In terms of the
periodic families of stable homotopy theory, if the homotopy of $B$
is $v_n$-periodic but not $v_{n+1}$-periodic, then frequently $K(B)$ is
$v_{n+1}$-periodic but not $v_{n+2}$-periodic.

Since the (fundamental) period $|v_{n+1}| = 2p^{n+1}-2$ of
$v_{n+1}$-periodicity is longer than the period $|v_n| = 2p^n-2$
of $v_n$-periodicity, we think of this phenomenon as an increase,
or lengthening, of wavelengths.  This is what we informally call a
``redshift''.  In a related fashion, the $v_{n+1}$-periodic phenomena
are usually hidden or nested behind the $v_n$-periodic ones, hence more
subtle and difficult to detect.  Again this corresponds informally to
less energetic light, propagating at lower frequencies.

The height filtration is also related to the sequence of Hopf subalgebras
$$
0 \subset \dots \subset \sE(n) = E(Q_0, \dots, Q_n) \subset \dots
$$
of the Steenrod algebra $\sA$, and their annihilating subalgebras
$$
\sA_* \supset \dots \supset
  (\sA//\sE(n))^* = P(\bar\xi_k \mid k\ge1) \otimes E(\bar\tau_k \mid
	k \ge n+1) \supset \dots \,.
$$
The latter nested sequence of $\sA_*$-comodule subalgebras are invariant
under the Dyer--Lashof operations that arise from thinking of the dual
Steenrod algebra $\sA_*$ as $H_*(H)$, where $H = H\bF_p$ is a commutative
structured ring spectrum.

\section{Redshift in the $K$-theory of rings}

We start with examples of chromatic redshift in the algebraic $K$-theory
of discrete rings.

Let $k$ be a finite field of characteristic~$p$, with algebraic closure
$\bar k$.  Quillen proved \cite{Q72}*{\S11} that $H_i(BGL(\bar k); \bF_p)
= 0$ for $i>0$, so that $K(\bar k)_p \simeq H\bZ_p$.  Furthermore, he
deduced that $\pi_* K(k)_p \cong \pi_* K(\bar k)_p^{hG_k}$ for $*\ge0$,
where the absolute Galois group $G_k$ acts continuously on $K(\bar k)$,
so $K(k)_p \simeq H\bZ_p$.  Multiplication by $p$ acts injectively on
$\pi_* K(\bar k)_p$, hence also on $\pi_* K(k)_p$.  Think of $p$ as a
lift of $p = v_0 \in \pi_* BP$, where $BP$ is the Brown--Peterson ring
spectrum with $\pi_* BP = \bZ_{(p)}[v_n \mid n\ge1]$.

For a separably closed field $\bar F$ of characteristic $\ne p$
(including $0$), Lichtenbaum conjectured that $\pi_t K(\bar F)_p$ is
$\bZ_p$ for $t\ge0$ even and $0$ for $t$ odd.  This was proved by Suslin
\cite{S83}*{Cor.~3.13}, and implies that $K(\bar F)_p \simeq ku_p$ and
$\hat L_1 K(\bar F) \simeq KU_p$.  Here $ku$ is the connective cover of
the complex topological $K$-theory ring spectrum $KU$, and $\hat L_n =
L_{K(n)}$ denotes Bousfield localization \cite{B79} with respect to the
Morava $K$-theory ring spectrum $K(n)$.  Multiplication by the Bott
element $u \in \pi_2 ku_p$ acts bijectively on $\pi_* K(\bar F)_p$,
for $*\ge0$.

Let $F$ be a number field, with a ring of $S$-integers $A$.
$$
\xymatrix{
A \ar[rr] && F \\
\bZ \ar[r] \ar[u] & \bZ[1/p] \ar[r] & \bQ \ar[u]
}
$$
Quillen conjectured \cite{Q75}*{\S9} that there is a spectral sequence
$$
E^2_{s,t} = H^{-s}_{\et}(\Spec A; \bZ_p(t/2))
	\Longrightarrow \pi_{s+t} K(A)_p
$$
converging in total degrees $\ge1$.  Here $H^*_{\et}(-)$ denotes
{\'e}tale cohomology, which is only well-behaved if $1/p \in A$, and
$\bZ_p(t/2) \cong \pi_t K(\bar F)_p$.  For $A = F$ this means that $\pi_*
K(F)_p \cong \pi_* K(\bar F)^{hG_F}_p$ for $*\ge1$, where $G_F$ is the
absolute Galois group.  The general case requires the more elaborate
framework of {\'e}tale homotopy types.  Passing to mod~$p$ homotopy,
a lift $\beta \in \pi_{2p-2}(S/p)$ of $u^{p-1} \in \pi_{2p-2}(ku;
\bZ/p)$ would act bijectively on $\pi_*(K(A); \bZ/p)$, for $*\ge1$.
Think of $\beta = v_1$ as a lift of $v_1 \in \pi_*(BP; \bZ/p)$.

Thomason \cite{T85}*{Thm.~4.1} proved Quillen's conjecture, up to the
localization given by inverting~$\beta$.  In particular, $\pi_*(K(F);
\bZ/p)[1/\beta] \cong \pi_*(K(\bar F)^{hG_F}; \bZ/p)$ for $*\ge2$.  It
remained to show that $\pi_*(K(A); \bZ/p) \to \pi_*(K(A); \bZ/p)[1/\beta]$
is an isomorphism for $*\ge2$.  Waldhausen \cite{W84}*{p.~193} noted that
this amounts to asking that $K(A) \to L_1 K(A)$ is a $p$-adic equivalence,
in sufficiently high degrees.  Here $L_n = L_{E(n)}$ denotes Bousfield
localization with respect to the Johnson--Wilson ring spectrum $E(n)$,
or equivalently with respect to $BP[1/v_n]$.

Using topological cyclic homology, Hesselholt--Madsen \cite{HM03}*{Thm.~A}
confirmed Quillen's conjecture for valuation rings in local number fields,
after special cases were treated by B{\"o}kstedt--Madsen \cite{BM94},
\cite{BM95} and Rognes \cite{R99a}, \cite{R99b}.

Finally, Voevodsky's proof \cite{V03}, \cite{V11} of the Milnor and
Bloch--Kato conjectures confirmed Quillen's conjecture for rings of
integers in global number fields.

\section{Redshift in the $K$-theory of ring spectra}

We continue with examples of chromatic redshift in the context
of algebraic $K$-theory of structured ring spectra.

Let $L = E(1)$ be the Adams summand of $KU_{(p)}$, and $\ell =
BP\langle1\rangle$ its connective cover.  Using topological cyclic
homology, Ausoni--Rognes \cite{AR02}*{Thm.~0.4} computed $V(1)_*
K(\ell_p)$, and Ausoni \cite{A10}*{Thm.~1.1} computed $V(1)_* K(ku_p)$,
where $p\ge5$ and $V(1) = S/(p,v_1)$ is the Smith-Toda spectrum of
chromatic type~$2$.  Using a localization sequence of Blumberg--Mandell
\cite{BM08}*{p.~157}, this also calculates $V(1)_* K(L_p)$ and $V(1)_*
K(KU_p)$.  In each case, a lift $v_2 \in \pi_{2p^2-2} V(1)$ of $v_2 \in
V(1)_* BP$ acts bijectively on the answer $V(1)_* K(B)$, for $*\ge 2p-2$.

The results are compatible with the existence of a spectral sequence
$$
E^2_{s,t} = H^{-s}_{\mot}(\Spec B; \bF_{p^2}(t/2))
	\Longrightarrow V(1)_{s+t} K(B)
$$
for suitable ``$\ell_p$-algebras of $S$-integers'' $B$, converging
in sufficiently high total degrees.  Here $H^*_{\mot}(-)$ refers to a
hypothetical form of motivic cohomology for commutative structured
ring spectra, and $\bF_{p^2}(t/2) \cong V(1)_t E_2$ where $E_2$
is the $K(2)$-local Lubin--Tate ring spectrum with $\pi_* E_2 =
\bW\bF_{p^2}[[u_1]][u]$.

The appearance of the field $\bF_{p^2}$
is needed to account for the sign in Ausoni's relation $b^{p-1} = -v_2$
in $V(1)_* K(ku_p)$, since if $b$ is represented by $\alpha u^{p+1}$
and $v_2$ by $u^{p^2-1}$ then $\alpha^{p-1} = -1$, which cannot be
satisfied for $\alpha \in \bF_p$.

\section{An analogue of the Lichtenbaum--Quillen conjectures}

Consider a Galois extension $L_p[1/p] \to M$, like in \cite{R08}*{\S4}.
By an $\ell_p$-algebra of integers in $M$ we mean a connected (only
trivial idempotents) commutative $\ell_p$-algebra $B$, with a structure
map to $M$, such that $B$ is semi-finite (retract of a finite cell
module), or perhaps dualizable, as an $\ell_p$-module:
$$
\xymatrix{
\Omega_1 \\
B \ar[rr] \ar[u] && M \\
\ell_p \ar[u] \ar[r] & L_p \ar[r] & L_p[1/p] \ar[u]_G \\
& J_p \ar[u]
}
$$
For $S$-integers we may allow localizations that invert $p$ or~$v_1$.
Let $\Omega_1$ be the $p$-completed homotopy colimit of all such $B$,
i.e., the $\ell_p$-algebraic integers.
% ((Does this exist?))

By analogy with Quillen's conjecture/Voevodsky's theorem we predict that
$v_2$ acts bijectively on $V(1)_* K(B)$, for $* \gg 0$.
By analogy with Lichtenbaum's conjecture/Suslin's theorem, we
predict that $V(1)_* K(\Omega_1) \cong V(1)_* E_2$, in all sufficiently
high degrees, and that $\hat L_2 K(\Omega_1) \simeq E_2$.

In the case when $B \to \Omega_1$ is an unramified $G$-Galois
extension, the hypothetical motivic cohomology would reduce
to group cohomology, and $V(1)_* K(B) \cong V(1)_* K(\Omega_1)^{hG}$
for $* \gg 0$.  The general case would require a more
elaborate construction than the familiar homotopy fixed points.
Even establishing the existence of a ring spectrum map $K(ku) \to E_2$
seems to be an open problem.

Similarly, for $n\ge1$ let $E_n$ be the $K(n)$-local Lubin--Tate ring
spectrum, and let $e_n$ be its connective cover, so that
$E_n = e_n[1/u]$.  Consider Galois
extensions $E_n[1/p] \to M$ and connected commutative $e_n$-algebras
$B$, with a structure map to $M$, such that $B$ is semi-finite as an
$e_n$-module:
$$
\xymatrix{
\Omega_n \\
B \ar[rr] \ar[u] && M \\
e_n \ar[u] \ar[r] & E_n \ar[r] & E_n[1/p] \ar[u] \\
& \hat L_n S \ar[u]
}
$$
Let $\Omega_n$ be the $p$-completed homotopy colimit of all such
$B$, i.e., the $e_n$-algebraic integers.

Let $F$ be a finite $p$-local spectrum admitting a $v_{n+1}$ self
map $v \: \Sigma^d F \to F$, cf.~Hopkins--Smith \cite{HS98}*{Def.~8}.
The finite localization functor $L^f_{n+1}$, which annihilates all finite
$E(n+1)$-acyclic spectra \cite{M92}*{Thm.~4}, is a smashing localization
such that $F_* L^f_{n+1} X \cong F_* X[1/v]$ for all spectra $X$.

I stated something like the following at Schlo{\ss} Ringberg in January
1999 and in Oberwolfach in September 2000:

\begin{conjecture} \label{conj}
Let $B \to \Omega_n$ and $(F, v)$ be as above.

(a)
Multiplication by $v$ acts bijectively on $F_* K(B)$ for $* \gg 0$,
and $K(B) \to L^f_{n+1} K(B)$ is a $p$-adic equivalence in
sufficiently high degrees.

(b)
There are isomorphisms $F_* K(\Omega_n) \cong F_* E_{n+1}$ for $* \gg 0$,
and $\hat L_{n+1} K(\Omega_n) \simeq E_{n+1}$.
\end{conjecture}

The cases $n=-1$ and $n=0$ correspond to Quillen's results and the
proven Lichtenbaum--Quillen conjectures, respectively.

\section{The cyclotomic trace map}

We can detect chromatic redshift in algebraic $K$-theory using
the cyclotomic trace map to topological cyclic homology, or one
of its variants.

The topological Hochschild homology $\THH(B)$ of a commutative $S$-algebra
$B$ is an $S^1$-equivariant spectrum whose underlying spectrum with
$S^1$-action can be constructed as $B \otimes S^1$, where $\otimes$
refers to the tensored structure of commutative $S$-algebras over spaces.
Let
$$
\THH(B)^{hS^1} = F(ES^1_+, \THH(B))^{S^1}
$$
be the $S^1$-homotopy fixed points of $\THH(B)$, and let
$$
\THH(B)^{tS^1} = [\widetilde{ES^1} \wedge F(ES^1_+, \THH(B))]^{S^1}
$$
be its $S^1$-Tate construction, also denoted $t_{S^1}
\THH(B)^{S^1}$ or $\hat\bH(S^1, \THH(B))$.  Here $ES^1$ is a free
contractible $S^1$-space, and $\widetilde{ES^1}$ is the mapping cone of
the collapse map $ES^1_+ \to S^0$.  Homotopy fixed point spectra model
group cohomology, and the Tate construction models Tate cohomology.

Think of $B$ as a ring spectrum of functions on a brave new scheme $X$.
Then $B \wedge \dots \wedge B$ is the ring of functions on $X \times
\dots \times X$, so $\THH(B)$ plays the role of the ring of functions on
the free loop space $\Map(S^1, X) = \Lambda X$, and $\THH(B)^{hS^1}$ is
like the ring of functions on the Borel construction $ES^1_+ \wedge_{S^1}
\Lambda X$.  The Tate construction is a periodicized version of the
Borel construction.

There is a natural trace map
$$
K(B) \longto \THH(B)
$$
that factors through the fixed point spectra $\THH(B)^{C_r}$ for
all finite subgroups $C_r \subset S^1$.  In particular, there is a
limiting map
$$
K(B) \longto TF(B; p) = \holim_n \THH(B)^{C_{p^n}} \,.
$$
Continuing with the canonical map from fixed points to homotopy fixed
points, the target of
$$
\holim_n \THH(B)^{C_{p^n}} \longto \holim_n \THH(B)^{hC_{p^n}}
$$
is $p$-adically equivalent to $\THH(B)^{hS^1}$.  The cyclotomic
structure of $\THH(B)$ gives a similar map
$$
\holim_n \THH(B)^{C_{p^n}} \longto \holim_n \THH(B)^{tC_{p^{n+1}}}
$$
whose target is $p$-adically equivalent to $\THH(B)^{tS^1}$.

The topological Hochschild construction itself does not introduce
a redshift, since $\THH(B)$ is a commutative $B$-algebra.  However,
in all the computations made so far, any $v_{n+1}$-periodicity that is
seen in the algebraic $K$-theory $K(B)$ has already been visible in the
$S^1$-Tate construction $\THH(B)^{tS^1}$.

Furthermore, it is possible to see in homological terms where
the redshift arises, in terms of these $S^1$-equivariant constructions.

\section{Circle-equivariant redshift}

The algebra $H_*(e_n)$ appears to be unwieldy for $n\ge2$, but
there is a map $BP\langle n\rangle \to e_n$ of (not necessarily
commutative) $S$-algebras, covering the usual map $E(n) \to E_n$,
and the augmentation $BP\langle n\rangle \to H$ induces 
an identification
$$
H_*(BP\langle n\rangle) \cong
	P(\bar\xi_k \mid k\ge1) \otimes E(\bar\tau_k \mid k \ge n+1)
$$
of subalgebras of the dual Steenrod algebra
$$
\sA_* = P(\bar\xi_k \mid k\ge1) \otimes E(\bar\tau_k \mid k\ge0) \,.
$$
Forgetting some structure, we can therefore think of the homology $H_*(B)$
of a commutative $e_n$-algebra~$B$ as a commutative $H_*(BP\langle
n\rangle)$-algebra.  This makes the Adams spectral sequence
$$
E_2^{s,t}(B) = \Ext_{\sA_*}^{s,t}(\bF_p, H_*(B))
	\Longrightarrow \pi_{t-s}(B^\wedge_p)
$$
an algebra over the Adams spectral sequence
$$
E_2^{s,t}
	= \Ext_{\sA_*}^{s,t}(\bF_p, H_*(BP\langle n\rangle))
	\Longrightarrow \pi_{t-s}(BP\langle n\rangle^\wedge_p)
$$
which collapses at the $E_2$-term
$$
E_2^{*,*} = P(v_0, \dots, v_n)
$$
and converges to the homotopy
$$
\pi_* BP\langle n\rangle^\wedge_p
	\cong \bZ_p[v_1, \dots, v_n] \,.
$$

The B{\"o}kstedt spectral sequence
$$
E^2_{s,t}(B) = HH_s(H_*(B))_t
	\Longrightarrow H_{s+t}(\THH(B))
$$
is then an algebra spectral sequence over
$$
E^2_{*,*} = HH_*(H_*(BP\langle n\rangle)) \cong
  H_*(BP\langle n\rangle) \otimes E(\sigma\bar\xi_k \mid k\ge1)
  \otimes \Gamma(\sigma\bar\tau_k \mid k \ge n+1)
$$
converging to $H_*(\THH(BP\langle n\rangle))$.  Here $\sigma$ denotes
the suspension operator, coming from the $S^1$-action on $\THH$, and
$\Gamma(x) = \bF_p\{\gamma_j x \mid j\ge0\}$ denotes the divided power
algebra on $x$.

The Dyer--Lashof operations $Q^{p^k}(\bar\tau_k) = \bar\tau_{k+1}$
in $\sA_*$ (coming from the commutative $S$-algebra structure
on $H$), imply multiplicative extensions $(\sigma\bar\tau_k)^p =
\sigma\bar\tau_{k+1}$, for $k \ge n+1$, which in turn imply that the Bockstein
images $\beta(\sigma\bar\tau_{k+1}) = \sigma\bar\xi_{k+1}$ vanish in
the abutment.  This argument, see Ausoni \cite{A05}*{Lem.~5.3}, implies
differentials
$$
d^{p-1}(\gamma_j \sigma\bar\tau_k)
	\doteq \sigma\bar\xi_{k+1} \cdot \gamma_{j-p} \sigma\bar\tau_k
$$
for all $j\ge p$, which leave
$$
E^p_{*,*} = E^\infty_{*,*} \cong H_*(BP\langle n\rangle)
	\otimes E(\sigma\bar\xi_1, \dots, \sigma\bar\xi_{n+1})
	\otimes P_p(\sigma\bar\tau_k \mid k \ge n+1)
$$
converging to
$$
H_*(\THH(BP\langle n\rangle)) \cong H_*(BP\langle n\rangle)
        \otimes E(\sigma\bar\xi_1, \dots, \sigma\bar\xi_{n+1})
	\otimes P(\sigma\bar\tau_{n+1}) \,.
$$
This will still have trivial $v_{n+1}$-periodic homotopy, but
note how building in a circle action gives rise to the class
$\sigma\bar\tau_{n+1}$.

The homological Tate spectral sequence
$$
E^2_{s,t}(B) = \hat H^{-s}(S^1; H_t(\THH(B)))
	\Longrightarrow H^c_{s+t}(\THH(B)^{tS^1})
$$
converges to a limit that we call the continuous homology of
$\THH(B)^{tS^1}$.  It is an algebra spectral sequence over
$$
E^2_{*,*} = \hat H^{-*}(S^1; H_*(\THH(BP\langle n\rangle)))
\cong P(t^{\pm1}) \otimes H_*(\THH(BP\langle n\rangle))
$$
converging to $H^c_*(\THH(BP\langle n\rangle)^{tS^1})$.
Here
$$
d^2(t^i \cdot x) = t^{i+1} \cdot \sigma x
$$
for all $x$, which leaves
\begin{multline*}
E^3_{*,*} = P(t^{\pm1}) \otimes 
P(\bar\xi_1^p, \dots, \bar\xi_{n+1}^p, \bar\xi_k \mid k\ge n+2) \\
\otimes E(\tau'_k \mid k\ge n+2)
\otimes 
E(\bar\xi_1^{p-1} \sigma\bar\xi_1, \dots,
	\bar\xi_{n+1}^{p-1} \sigma\bar\xi_{n+1})
\end{multline*}
where $\tau'_k = \bar\tau_k - \bar\tau_{k-1}(\sigma\bar\tau_{k-1})^{p-1}$
for each $k\ge n+2$.  Note that $\bar\tau_{n+1}$ supports a
nontrivial $d^2$-differential to $t \cdot \sigma\bar\tau_{n+1}$, and
does not survive to the $E^\infty$-term, while the $\tau'_k$ for
$k\ge n+2$ are $d^2$-cycles, due to the known multiplicative extension.

This spectral sequence often collapses at this stage
\cite{BR05}*{Prop.~6.1}, and there can be $\sA_*$-comodule extensions
that combine $p^{n+1}$ shifted copies of
$$
P(\bar\xi_1^p, \dots, \bar\xi_{n+1}^p, \bar\xi_k \mid k\ge n+2)
\otimes E(\tau'_k \mid k\ge n+2)
$$
to a copy of $P(\bar\xi_k \mid k\ge 1) \otimes E(\tau'_k \mid k\ge
n+2) \cong H_*(BP\langle n+1\rangle)$.  The PhD theses of Sverre
Lun{\o}e--Nielsen \cite{LR12}, \cite{LR11} and Knut Berg (to appear)
address these questions.  Note the transition from $H_*(BP\langle
n\rangle)$ to $H_*(BP\langle n+1\rangle)$, with non-trivial
$v_{n+1}$-periodic homotopy groups.
The typical result is that $H^c_*(\THH(B)^{tS^1})$ is an
algebra over $H^c_*(\THH(BP\langle n\rangle)^{tS^1})$, which
has an associated graded of the form
$$
P(t^{\pm p^{n+1}}) \otimes H_*(BP\langle n+1\rangle)
	\otimes E(\nu_1, \dots, \nu_{n+1})
$$
where $\nu_k$ is a $t$-power multiple of $\bar\xi_k^{p-1} \sigma\bar\xi_k$,
but that there is room for further $\sA_*$-comodule extensions.

This implies that the inverse limit Adams spectral sequence
$$
E_2^{s,t}(B) = \Ext_{\sA_*}^{s,t}(\bF_p, H^c_*(\THH(B)^{tS^1}))
	\Longrightarrow \pi_{t-s} \THH(B)^{tS^1}_p
$$
is an algebra over the Adams spectral sequence
$$
E_2^{s,t} = \Ext_{\sA_*}^{s,t}(\bF_p, H^c_*(\THH(BP\langle n\rangle)^{tS^1}))
	\Longrightarrow \pi_{t-s} \THH(BP\langle n\rangle)^{tS^1}_p
$$
which contains factors like
$$
\Ext_{\sA_*}^{*,*}(\bF_p, H_*(BP\langle n+1\rangle))
	\cong P(v_0, \dots, v_n, v_{n+1}) \,.
$$
Due to the exterior factors $E(\nu_1, \dots, \nu_{n+1})$ there is room
for differentials that might truncate the periodic $v_{n+1}$-action
visible above, but empirically this does not happen.  A theory that
explains the general picture is, however, currently lacking.

\section{Beyond elliptic cohomology}

Do $K(tmf)$ and $\THH(tmf)^{tS^1}$ detect $v_3$-periodic families?
Work in progress for $p=2$ with Bruner (2008).

\section{Waldhausen's localization tower}

The chromatic localization functors ($L_n$ and) $\hat L_n$ and the finite
localizations functors $L^f_n$ fit in a diagram of commutative structured
ring spectra
$$
\xymatrix{
& & E_n & & & KU_p \\
& & \hat L_n S \ar[u]^-{\bG_n} & & & J_p \ar[u]^{\bZ_p^\times}
	& H\bQ \\
S_{(p)} \ar[r] & \dots \ar[r] & L^f_n S \ar[r] \ar[u]
	& L^f_{n-1} S \ar[r] & \dots \ar[r]
	& L^f_1 S \ar[r] \ar[u]
	& L^f_0 S \ar[u]_{\simeq}
}
$$
where $L^f_n S \to L_n S$ is an equivalence for $n\le1$, but probably
not for $n\ge2$, according to the wisdom concerning Ravenel's telescope
conjecture~\cite{MRS01}.  Applying algebraic $K$-theory
to the lower row one gets a telescopic localization tower
$$
K(S_{(p)}) \longto \dots \longto K(L^f_n S) \longto K(L^f_{n-1} S)
	\longto \dots \longto K(L_1 S) \longto K(\bQ)
$$
similar to that of \cite{W84}*{p.~174}, interpolating between the
geometrically significant algebraic $K$-theory of spaces on the left
hand side, and the arithmetically significant algebraic $K$-theory of
number fields on the right hand side.  Waldhausen worked with $L_n$,
and explicitly assumed that it is a finite localization functor, but
we can work with $L^f_n$ instead.  This ensures that each finite cell
$L^f_n S$-module is $L^f_n$-equivalent to a finite cell $S$-module,
as can be proved by induction on the number of $L^f_n S$-cells.

Let $\sC_0$ be the category of finite $p$-local spectra, and let $w_n
\sC_0$ be the subcategory of $E(n)_*$-equivalences, or equivalently
of $L^f_n$-equivalences, for $n\ge0$.  Let $\sC_n = \sC_0^{w_{n-1}}$
denote the full subcategory of $E(n-1)_*$-acyclic spectra, i.e., the
finite spectra of type $\ge n$, for $n\ge1$.  Then $K(\sC_0, w_n) \simeq
K(L^f_n S)$, and Waldhausen's localization theorem \cite{W84}*{\S3}
recognizes the homotopy fiber of $K(L^f_n S) \to K(L^f_{n-1} S)$ as
$K(\sC_n, w_n)$, i.e., the algebraic $K$-theory of finite spectra of
type $\ge n$, with respect to the $E(n)_*$-equivalences.  We get a
homotopy fiber sequence
$$
K(\sC_n, w_n) \longto K(L^f_n S) \longto K(L^f_{n-1} S) \,.
$$

Let $\sK_n^{sm}$ be the category of small $K(n)$-local spectra, and
let $\sK_n'$ be the full subcategory of $K(n)$-localizations of finite
spectra of type $\ge n$.  Hovey--Strickland \cite{HS99}*{Thm.~8.5}
show that the inclusion $\sK_n' \subset \sK_n^{sm}$ is an idempotent
completion, so the induced map $K(\sK_n') \to K(\sK_n^{sm})$ induces an
isomorphism on $\pi_i$ for each $i\ge1$.  The localization functors $L_n$
and $\hat L_n$ agree on $\sC_n$, hence induce an equivalence $K(\sC_n,
w_n) \simeq K(\sK_n')$.  Thus we have a map
$$
K(\sC_n, w_n) \longto K(\sK_n^{sm}) \,,
$$
which induces a $\pi_i$-isomorphism for each $i\ge1$.  We view
$\sK^{sm}_n$ as a category of suitably small $\hat L_n S$-modules.
% ((What happens at the level of $\pi_0$?))

Let $\sE_n^{df}$ be the category of $E_n$-module spectra that have
degreewise finite homotopy groups.  Base change along the $K(n)$-local
pro-$\bG_n$-Galois extension $\hat L_n S \to E_n$ takes $\sK_n^{sm}$
to $\sE_n^{df}$, and conversely \cite{HS99}*{Cor.~12.16}, so it is
plausible that a Galois descent comparison map
$$
K(\sK_n^{sm}) \longto K(\sE_n^{df})^{h\bG_n}
$$
is close to an equivalence.  Finally, $K(\sE_n^{df})$ is related
to the algebraic $K$-theory of $E_n$ and its various localizations.
For $n=1$ we have $E_1 = KU_p$, and $K(\sE_1^{df})$ is the algebraic
$K$-theory of $p$-nilpotent finite cell $KU_p$-modules, which sits
\cite{B?}*{Prop.~11.15} in a homotopy fiber sequence
$$
K(\sE_1^{df}) \longto K(KU_p) \longto K(KU_p[1/p]) \,.
$$
In general, this fiber sequence is replaced by an $n$-dimensional
cubical diagram.  Note that the transfer map $K(KU/p) \to K(\sE_1^{df})$
associated to $KU_p \to KU/p$ is far from an equivalence, by the
calculations of \cite{AR12}*{Cor.~1.3}, so there does not appear to be
any easy way to describe the algebraic $K$-theory of degreewise finite
$E_n$-modules in terms of d{\'e}vissage, cf.~\cite{W84}*{p.~188}.

$$
\xymatrix{
& K(\sE_n^{df}) \ar[r] & K(E_n) \\
K(\sC_n, w_n) \ar[dr] \ar[r] & K(\sK_n^{sm}) \ar[u]^-{\bG_n} \\
\dots \ar[r] & K(L^f_n S) \ar[r] & K(L^f_{n-1} S) \ar[r] & \dots
}
$$
Conjecture~\ref{conj} about the structure of the algebraic $K$-theory of
$E_n$ (and various localizations) is therefore also a statement about
$K(\sE_n^{df})$, and conjecturally about $K(\sK_n^{sm})$, which rather
precisely measures the difference between $K(L^f_n S)$ and $K(L^f_{n-1}
S)$.

\section{The spherical case}

Calculations of $TC(S; p)$, $K(\bZ)$ and $TC(\bZ; p)$ were assembled to
a calculation of $K(S)$ at $p=2$ in \cite{R02} and at odd regular
primes in \cite{R03}.  These results describe the cohomology of $K(S)$
as an $\sA$-module in all degrees (up to an extension in the odd case),
and lead to Adams spectral sequence calculations in a finite range
of degrees.

The algebraic $K$-groups of $S$ are at least as complicated as those
of its stable homotopy groups.  The complex cobordism spectrum $MU$
has turned out to be a convenient halfway house
$$
S \longto MU \longto H
$$
between homology and homotopy.  The Thom equivalence $MU \wedge MU
\simeq MU \wedge BU_+$ makes $S \to MU$ a Hopf--Galois extension
\cite{R08}*{\S12}, and
the cosimplicial Amitsur resolution
$$
[q] \longmapsto MU \wedge MU^{\wedge q}
$$
of $S$ is equivalent to the cobar resolution
$[q] \longmapsto MU \wedge BU^q_+$
for the $S[BU] = \Sigma^\infty(BU_+)$-comodule algebra $MU$.
Applying algebraic $K$-theory, an analogue of Quillen's conjecture
would predict that $K(S)$ is well approximated by the totalization
of the cosimplicial spectrum
$$
[q] \longmapsto K(MU \wedge MU^{\wedge q})
$$
rewriteable as $[q] \longmapsto K(MU \wedge BU^q_+)$.
If, by analogy with the Galois case, there are compatible maps $K(MU
\wedge BU^q_+) \to K(MU) \wedge BU^q_+$, then this might in turn be
approximated by the totalization of the cobar resolution
$[q] \longmapsto K(MU) \wedge BU^q_+$
for an $S[BU]$-comodule algebra structure on $K(MU)$.

Conceivably, this leads to a more conceptual understanding of $\pi_* K(S)$
in terms of $\pi_* K(MU)$ and Hopf--Galois descent, by analogy with the
Adams--Novikov spectral sequence description of $\pi_* S$ in
terms of $\pi_* MU$ and its $H_*(BU)$-coaction.  This has been a
motivating factor for the study of $K(MU)$, advertised in \cite{BR05}
and \cite{R09}, and pursued in \cite{LR11}.

\section{Higher redshift}

For a Lie group $G$ of rank~$k$, consider $(B \otimes G)^{hG}$ or
something like $(B \otimes G)^{tG}$.  If $B$ is $v_n$-periodic
but not $v_{n+1}$-periodic, then apparently $(B \otimes G)^{tG}$
is $v_{n+k}$-periodic.   Tested for $B = H$ and $G = T^k$ for
small~$k$, as well as for $G = SO(3)$ and $G = S^3$.
Work in progress (Rognes, 2008--2011) and in Torleif Veen's PhD
thesis (2013).

\end{document}